\documentclass[12pt, a4paper]{article}
\usepackage[francais]{babel}
\usepackage[latin1]{inputenc}
\usepackage{amssymb}
\usepackage{amsmath}
\usepackage{amsfonts}
\usepackage{mathrsfs}

\addto\captionsfrench{}
\addto\captionsfrench{}

\newcounter{chapitre}
\newtheorem{thm}{Theorem}[section]
\newtheorem{cor}[thm]{Corollary}

\newtheorem{prop}[thm]{Proposition}
\newtheorem{defn}{Definition}[section]

\newenvironment{dem}{{\noindent { $Proof.$ }}}

\newcommand{\N}{\mathbb{N}}

\def\keywords#1{{ \it{Keywords:}} {#1}} %

\begin{document}

\title{ A NEW FORMULA OF $q$-FUBINI NUMBERS VIA GON$\breve{C}$AROV POLYNOMIALS  }


%
%
\author{Adel Hamdi }

\maketitle

\begin{center}

\centerline{\footnotesize{Faculty of Science of Gabes, Department of
Mathematics,}} \centerline{\footnotesize{ Cit\'e Erriadh 6072,  Zrig, Gabes,
Tunisia}}
\end{center}

\begin{abstract}
Connected the generalized Gon$\breve{c}$arov polynomials associated to a pair ($\partial,\mathcal{Z}$)
of a delta operator $\partial$ and an interpolation grid $\mathcal{Z}$, introduced by Lorentz, Tringali and Yan in \cite{2}, with the theory of binomial enumeration and order statistics, a new $q$-deformed of these polynomials given in this paper allows us to derive a new combinatorial formula of $q$-Fubini  numbers. A combinatorial proof and some nice algebraic and analytic properties have been expanded to the $q$-deformed version.
\end{abstract}
\begin{center}
.
\end{center}


\section{Introduction}
This paper grew out of the recent work, generalized Gon$\breve{c}$arov polynomials,   of Lorentz, Tringali and Yan in \cite{2} where these polynomials are seen as a basis of solutions for the Interpolation problem: Find a polynomial $f$($x$) of degree $n$ such that the $i$th delta operator $\partial$ of $f$($x$) at a given complex number $a_i$ has value $b_i$, for $i=0,1,2,...$. There is a natural $q$-analog of this interpolation by replacing the delta operator with a $q$-delta operator,  we extend these polynomials into a generalized $q$-Gon$\breve{c}$arov basis  $(t_{n,q}(x))_{n\geq0}$, defined by the $q$-biorthogonality relation $\varepsilon_{z_i}(\partial_q^i(t_{n,q}(x)))=[n]_q!\delta_{i,n}$, for all $i,n\in \mathbb{N}$, where $\mathcal{Z} =$ ($z_i$)$_{i\geq0}$ is a sequence of scalars and $\varepsilon_{z_i}$ the evaluation at $z_i$. Using ordered partitions, we give a combinatorial formula and a combinatorial proof of the constant terms of $t_{n,q}(x)$ where there are inspired from that of $t_n$($x$) in \cite{2}. This motivation us to allow to derive a new combinatorial formula of the $q$-Fubini numbers.

In particular, we described the algebraic properties of the sequence of
$q$-Gon$\breve{c}$arov polynomials which are based on the observations in \cite{1,2,3}.
A sequence of polynomials ($f_n(x)$)$_{n\geq0}$ $q$-biorthogonal  to a sequence of $q$-shift-invariant
operators ($\Psi$)$_{i\geq0}$ of the form $\sum_{j\geq 0} \lambda_j^{(i)}\partial_q^{i+j}$, where
$\lambda_j^{(i)}\in \mathbb{K}$ and $\lambda_0^{(i)}\neq 0$, is the unique sequence
satisfies the following relation:
\begin{equation}
\varepsilon_{0}\Psi_i(f_n(x))=[n]_q!\ \delta_{i,n},\qquad \mbox {for each}\ i=0,1,...,n.
\end{equation}

The theory of polynomials of $q$-binomial type plays a fundamental role in $q$-umbral
calculus, or finite
operator calculus pioneered by M. E. H. Ismail, \cite{6}. Each polynomial sequence of
$q$-binomial type can be characterized by a linear operator called $q$-delta operator and a
basic principles of $q$-delta operator exhibits fascinating relationships between the
formal power series methods of combinatorics.

 The main objective of this paper  is to present the algebraic  properties of generalized
 $q$-Gon$\breve{c}$arov polynomials, we will give a combinatorial formula for its constant terms and we derive a new combinatorial formula of $q$-Fubini numbers. The rest of
 the paper is organized as follows. We begin, in Section $2$, by giving an outline of the
 theory of sequences of polynomials $q$-biorthogonal to a sequence of linear functionals.
 So, this section  contains the definition and basic properties of a $q$-delta operators
 and their associated sequence of basic polynomials, the details can be found in \cite{6}.
 In Section $3$, we introduce the sequence of generalized $q$-Gon$\breve{c}$arov
 polynomials $t_{n,q}(x;\partial_q,\mathcal{Z})$ associated with a $q$-delta operator
 $\partial_q$ and a grid $\mathcal{Z}$ and we will describe its algebraic properties and
 characterizations.
In Section 4, we present a combinatorial description of the constant terms of the generalized
$q$-Gon$\breve{c}$arov polynomials in terms of rankings on ordered partitions with some
conditions. We conclude this paper by a new combinatorial formula of the $q$-Fubini numbers.

\section{$q$-delta operator and basic polynomials}

In this section, we begin by recalling some $q$-notions of quantum analysis adapting the
notation considered in \cite{9}.\\
For any complex number $a$, the basic number and $q$-factorial are defined as
$$[a]_q := \frac{1-q^a}{1-q} , \;\  q \neq 1; \;\ [n]_q!:=[n]_q [n - 1]_q...[2]_q [1]_q ,\;\ n = 1, 2,...,$$
with $[0]_q ! := 1$ and the scalar $q$-shifted factorials are defined as
$$(a, q)_0 := 1, \;\ (a, q)_n :=\prod^{n-1}_{k=0}(1 - aq^k ),\;\  n = 1, 2,..., $$
also we define
$${n \brack k}_q:=\frac{[n]_{q}!}{[k]_{q}![n-k]_{q}!},\;\ 0\leq k \leq n.
 $$
\vskip0.5cm
Since the generalized Gon$\breve{c}$arov polynomials are a sequence of binomial  type which can be characterized by a  linear operator called delta operator.
Then for to study the $q$-delta operator, Ismail  \cite{6} introduced an interesting $q$-analog of the  translation operator which is related to the continuous $q$-Hermite polynomials and extended it by
linearity to all polynomials.\\

The continuous $q$-Hermite polynomials are generated by
$$H_0(x|q)=1,\qquad H_1(x|q)= 2x,$$
$$H_{n+1}(x|q)= 2xH_n(x|q)-(1-q^n)H_{n-1}(x|q).$$
The action of the $q$-translation by $y$, $E_q^y$, on $H_n
(x|q)$ is defined by
$$E_q^y H_n(x|q)=H_n(x \oplus y|q):=\sum_{m\geq0} {n \brack m}_q H_m(x|q)
g_{n-m}(y)q^{(m^2-n^2)/4},$$
where the sequence of polynomials $(g_n(y))_{n\geq0}$ is
$$\frac{g_n(y)}{(q;q)_n}=\sum_{k=0}^{\lfloor n/2\rfloor}\frac{q^k}{(q^2;q^2)_k}
\frac{H_{n-2k}(y|q)}{(q;q)_{n-2k}}q^{(n-2k)^2/4}.$$
Now, we recall the definitions of a $q$-shift-invariant operator and a $q$-delta operator.
\begin{defn}
$(i)$\ A $q$-shift-invariant operator $\Lambda_q$ is a linear mapping polynomials to
polynomials which commutes with $E_q^a$
 for all complex numbers $a$.\\
$(ii)$\  A $q$-delta operator $\partial_q$
is a $q$-shift-invariant operator satisfying $\partial_q x$ = a for some nonzero constant a.
\end{defn}

By the above definition, $\partial_q$ has many properties such that deg($\partial_q$($f$))$ =$ deg($f$) $-1$ for any $f \in \mathbb{K}[x]$ and $\partial_q$($b$)$=0$ for every $q$-constant $b$.
\begin{defn}
Let $\partial_q$ be a $q$-delta operator. A polynomial sequence
$(p_n(x))_{n\geq0}$ is called the sequence of basic polynomials, or the associated basic sequence of
$\partial_q$ if
\begin{center}
\begin{description}
	\item[(i)] $p_0(x)=1$;
	\item[(ii)] degree of $p_n(x$) is n and $\tilde{g}_n(0)=0$, $n\geq1$, where
	\begin{equation}\label{eq1}
	\tilde{g}_n(x):=\sum_{k=0}^{\lfloor n/2\rfloor}\frac{(q;q)_n\ q^k}{(q^2;q^2)_k}
\frac{p_{n-2k}(x)}{(q;q)_{n-2k}};
	\end{equation}
  \item[(iii)] $\partial_q(p_n(x))=[n]_q p_{n-1}(x)$.
\end{description}
\end{center}
\end{defn}
\quad

So, from Theorem 14.7.1 of [2], every $q$-delta operator has a unique sequence of basic polynomial, which is a sequence of binomial type that satisfies
$$p_n(x \oplus y):=\sum_{m,j\geq0} {n \brack m}_q {n-m \brack 2j}_q (q;q^2)_j q^j p_m(y)
p_{n-m-2j}(x),$$
for all $n$.\\
Conversely, Theorem 14.7.2 of [2] implies that every sequence of polynomials of $q$-binomial type  is the basic
sequence for some $q$-delta operator.

 Let $\Lambda_q$ be a $q$-shift-invariant operator and $\partial_q$ a $q$-delta operator with  basic sequence of polynomials $(p_n(x))_{n\geq0}$. Then $\Lambda_q$ can be expanded as a formal power series of $\partial_q$, as
\begin{equation}
\Lambda_q=\sum_{k\geq0}\frac{a_k}{(q;q)_k} \partial_q^k,\qquad a_k:=
\varepsilon_0(\Lambda_q(\tilde{g}_k(x))).
\end{equation}
The polynomial $f$($t$) $= \sum_{k\geq0}\frac{a_k}{(q;q)_k}t^k$ is said the $\partial_q$-indicator of $\Lambda_q$ and we have the correspondence \begin{equation}
 f(t)=\sum_{k\geq0}\frac{a_k}{(q;q)_k} t^k \longleftrightarrow
\Lambda_q=\sum_{k\geq0}\frac{a_k}{(q;q)_k} \partial_q^k
\end{equation}
is an isomorphism from the ring $\mathbb{K}[\![t]\!]$ of formal power series into the ring of $q$-shift-invariant operators. So, a $q$-shift-invariant $\Lambda_q$ is invertible if and only if its
$\partial_q$-indicator $f(t)$ satisfies $f(0)\neq 0$  and it is a $q$-delta operator if and only if $f(0)=0$ and $f'(0)\neq 0$. For more details see [2].
\section{Generalized $q$-Gon$\breve{c}$arov Polynomials and Properties}

  In this  section, we give the main result which is the generalized $q$-Gon$\breve{c}$arov polynomials
  starting from the $q$-biorthogonality condition in the Gon$\breve{c}$arov interpolation
  problem. The classical version, means $q$ tends to $1$, was discussed in \cite{1} with
  the differentiation operator, in \cite{2} with an arbitrary delta operator and an
  explicit description for the $q$-difference operator is given in \cite{3}. Here we extend
  this theory to general $q$-delta operator.\\
	
 In the following, let $\partial_q$ be a $q$-delta operator with the basic sequence
 $(p_n(x))_{n\geq0}$ and let $(\Psi_i)_{i\geq0}$ be a sequence of $q$-shift-invariant
 operators where the $i^{th}$ term is defined by:
\begin{equation}
 \sum_{j\geq0} \lambda_j^{(i)} \partial_q^{i+j},
\end{equation}
where $\lambda_j^{(i)}\in \mathbb{K}$ and $\lambda_0^{(i)}\neq0$.

 A polynomial sequence $(p_n(x))_{n\geq0}$  is  called $q$-biorthogonal to the sequence of
 operators ($\Psi_i$)$_{n\geq0}$
if
\begin{equation}
\varepsilon_0(\Psi_i(p_n(x)))=[n]_q!\delta_{i,n}.
\end{equation}

\begin{thm}
For a sequence of  $q$-shift-invariants operators $(\Psi_i)_{n\geq0}$, there
exists a unique  polynomial sequence $(f_n(x))_{n\geq0}$, $f_n(x)$ is of degree n, that is $q$-biorthogonal to it.
Moreover, this sequence forms a basis of $\mathbb{K}[x]$: for every $f(x)\in \mathbb{K}[x]$
it holds
$$f(x)=\sum_{j=0}^{deg(f)}\frac{\varepsilon_0(\Psi_i(f(x)))}{[j]_q!} f_j(x).$$
\end{thm}

To prove this, follow the same technique adopted in \cite{1}, section 2, by replacing $n!$
with $[n]_q!$ (which tends to $n!$ as $q$ tends to $1$).\qquad\\

Let $\partial_q$ be a $q$-delta operator and $\mathcal{Z}=(z_i)_{i\geq0}$ be
a sequence of scalars that we called simply a grid and the values $z_i$ are
the nodes of $\mathcal{Z}$. For $a\in \mathbb{K}$, $E_q^a$ is an invertible
$q$-shift-invariant operator and from Theorem.14.7.4 and Theorem.14.7.5 of [2], there exist a unique $f_a(t)\in \mathbb{K}[\![t]\!]$ with
$f_a(0)\neq0$ such that $E_q^a=f_a(\partial_q)$.

\begin{defn}  The sequence of the generalized $q$-Gon$\breve{c}$arov polynomials
$(t_{n,q}(x)_{n\geq0}$ associated with the pair $(\partial_q, Z)$ is the sequence of polynomials $q$-biorthogonal to the sequence of
operators $(E_q^{z_i} \partial_q^i)_{i\geq0}$, means
\begin{equation}\label{eq2}
\varepsilon_{z_i}(\partial_q^i(t_{n,q}(x)))=[n]_q!\delta_{i,n}.
\end{equation}
Moreover, for any $f(x)\in \mathbb{K}[x]$, we have
\begin{equation}\label{base}
f(x)=\sum_{i=0}^{deg(f)}\frac{\varepsilon_{z_i}(\partial_q^i(f))}{[i]_q!} t_{i,q}(x).
\end{equation}
\end{defn}

The generalized $q$-Gon$\breve{c}$arov polynomials form a basis for the solutions of the
following interpolation
problem with a $q$-delta operator $\partial_q$, which we call the generalized
$q$-Gon$\breve{c}$arov
interpolation problem:\\

\textbf{Problem:} Given two sequences $z_0, z_1, ..., z_n$ and
$b_0, b_1, ... , b_n$ of real or complex numbers and a $q$-delta operator $\partial_q$,
find a (complex)
polynomial $p(x)$ of degree $n$ such that
\begin{equation}
\varepsilon_{z_i} \partial_q^i(p_n(x))=b_i,\qquad i=0,1,..., n.
\end{equation}

Let us emphasize that this polynomials depends only on $\partial_q$ and $\mathcal{Z}$, so
sometimes we can write $t_{n,q}(x,\partial_q,\mathcal{Z})$ instead of $t_{n,q}(x)$ and we
call it the generalized $q$-Gon$\breve{c}$arov basis associated with the pair
$(\partial_q,\mathcal{Z})$.

In the following result, we give some algebraic properties of generalized
$q$-Gon$\breve{c}$arov polynomials through which we can learn more about this polynomial.
\begin{prop}\label{pr1}
If $(t_{n,q}(x))_{n\geq0}$ is the generalized $q$-Gon$\breve{c}$arov basis associated with
the pair $(\partial_q,\mathcal{Z})$. Then $t_{0,q}(z_0)=1$ and $t_{n,q}(z_0)=0$ for all
$n\geq1$
\end{prop}
\begin{dem}
It's an immediate consequence of (\ref{eq2}). \hfill$\square$
\end{dem}
\vskip0.5cm

Let $\mathcal{Z}=(z_i)_{i\geq0}$ be a fixed grid. We will follow the same notation as in
\cite{2}, section 3, and we denote by $\mathcal{Z}^{(j)}$ and $\mathcal{O}$, respectively,
the grid whose $i$-th term is the $(i+j)$-th term  of $\mathcal{Z}$ (i.e.
$z_i^{(j)}=z_{i+j}$)  and the grid zero (i.e. $z_i=0$ for all $i$).
\vskip0.5cm

We give in the next proposition a $q$-analogue of the generalization of the differential
relations given in \cite{2}.
\begin{prop}\label{derv}
Let a fixed $j\in \mathbb{N}$ and for each $n\in \mathbb{N}$, we define the polynomial
$t_{n,q}^{(j)}$ as follows:
\begin{equation}\label{defG}
t_{n,q}^{(j)}(x):= \frac{(1-q)^j}{(q^{n+1};q)_{j}} \partial_q^j t_{n+j,q}(x),
\end{equation}
where $(t_{n,q}(x))_{n\geq0}$ is the generalized $q$-Gon$\breve{c}$arov basis associated
with the pair $(\partial_q,\mathcal{Z})$.\\ Then, $(t_{n,q}^{(j)}(x))_{n\geq0}$ is the
generalized $q$-Gon$\breve{c}$arov basis associated with the pair
$(\partial_q,\mathcal{Z}^{(j)})$.\\

In particular, we have
\begin{equation}\label{defg}
\partial_q^j t_{n,q}(x)=\frac{(q^{n-j+1},q)_{j}}{(1-q)^j} t_{n-j,q}^{(j)}(x).
\end{equation}
\end{prop}
\begin{dem} From Eq. (\ref{defG}), it is clear that $t_{n,q}^{(j)}(x)$ is a polynomial of
degree $n$.\\ Let $ i,n\in
\mathbb{N}$ with $0\leq i\leq n$. By the definition of the generalized $q$-Gon$\breve{c}$arov basis, it suffices to prove that
$$\varepsilon_{z_i^{(j)}}(\partial_q^i(t_{n,q}^{(j)}(x)))=[n]_q!\delta_{i,n},$$
where $z_i^{(j)}$ is the $i$-th node of the grid $\mathcal{Z}^{(j)}$. \\We have
\begin{align*}
\varepsilon_{z_i^{(j)}}(\partial_q^i(t_{n,q}^{(j)}(x)))&=\frac{(1-q)^j}{(q^{n+1};q)_{j}}
\varepsilon_{z_{i+j}}(\partial_q^{i+j}t_{n+j,q}(x))\\
&=\frac{(1-q)^j}{(q^{n+1};q)_{j}}[n+j]_q!\delta_{i+j,n+j}\\
&=\frac{(1-q)^j}{(q^{n+1};q)_{j}}\frac{(q;q)_{n+j}}{(1-q)^{n+j}}\delta_{i,n}\\
&= \frac{(q;q)_{n}}{(1-q)^n}\delta_{i,n}\\
&=[n]_q!\delta_{i,n},
\end{align*}
which follows from the identities $(a;q)_{n+m}=(a;q)_{n}(aq^{n};q)_{m}$,
$\delta_{i+j,n+j}=\delta_{i,n}$ and using the Eq. (\ref{eq2}). \hfill$\square$
\end{dem}

In the following, we investigate the relation between the  generalized $q$-Gon$\breve{c}$arov
polynomials and the basic polynomials of $q$-delta operator, present a  condition under
which $t_{n,q}(x,\partial_q,\mathcal{Z})$ is of $q$-binomial type.

\begin{prop}\label{basic}
The basic sequence of the $q$-delta operator $\partial_q$ is the generalized
$q$-Gon$\breve{c}$arov basis associated with the pair $(\partial_q,\mathcal{O})$.
\end{prop}
\begin{dem}
Let ($p_n(x)$)$_{n\geq0}$ be the basic sequence of the $q$-delta operator $\partial _q$. Then by the definition
$\partial_q(p_n(x))=[n]_q p_{n-1}(x)$. This implies that
$ \partial^i_q(p_n(x))= \frac{(q;q)_n}{(q;q)_{n-i}(1-q)^i} p_{n-i}(x).$
So, $\varepsilon_0(\partial^i_q(p_n(x)))= [n]_q!\delta_{i,n}$.$\hfill\square$
\end{dem}

\begin{cor}\label{equi}
A polynomial sequence $(p_n(x))_{n\geq0}$ is of $q$-binomial type if and only if it is the
generalized $q$-Gon$\breve{c}$arov basis associated with the pair
$(\partial_q,\mathcal{O})$ for a suitable choice of $\partial_q$.
\end{cor}

\begin{dem}
The direct sense implies that is a basis
sequence for fixed $q$-delta operator $\partial_q$. In view of the above proposition,
$(p_n(x))_{n\geq0}$ is the generalized $q$-Gon$\breve{c}$arov basis associated with the
pair $(\partial_q,\mathcal{O})$.

Conversely, Let $\partial_q$ be a $q$-delta operator and assume that ($p_n(x)$)$_{n\geq0}$
be the generalized $q$-Gon$\breve{c}$arov basis associated with the pair
$(\partial_q,\mathcal{O})$.
Based on Proposition \ref{derv}, we denotes by ($p_n^{(1)}(x)$)$_{n\geq0}$ the generalized
$q$-Gon$\breve{c}$arov basis associated with the pair $(\partial_q,\mathcal{O}^{(1)})$, which implies that, for all $n\geq0$,
$$\partial_q(p_n(x))=\frac{(q^n;q)_1}{1-q}p_{n-1}^{(1)}(x)
=[n]_q p_{n-1}(x),$$ since $\mathcal{O}^{(1)} = \mathcal{O}$.
 This together with Proposition \ref{pr1}  completes the proof.\hfill$\square$
\end{dem}
\vskip0.5cm

The next proposition gives the $q$-analogue of the extension of the shift-invariance
property studied in \cite{2}.
\begin{prop}\label{shiftt}
Let $(t_{n,q}(x))_{n\geq0}$ and $(h_{n,q}(x))_{n\geq0}$ be two  generalized
$q$-Gon$\breve{c}$arov bases associated with the pairs $(\partial_q,\mathcal{Z})$ and
$(\partial_q,\mathcal{W})$, respectively, with $\mathcal{W}=(w_i)_{i\geq0}$ be a
$q$-translation of the grid $\mathcal{Z}$ by $\xi\in \mathbb{K}$, i.e. $w_i=z_i \oplus
\xi$. Then
$h_{n,q}(x \oplus \xi)=t_{n,q}(x)$ for all $n\geq0$.
\end{prop}
\begin{dem}
By the uniqueness of the generalized $q$-Gon$\breve{c}$arov basis associated with the pair
$(\partial_q,\mathcal{Z})$ and since $h_{n,q}(x \oplus \xi)=E_q^{\xi}(h_{n,q}(x))$ is a
polynomial of degree $n$, it suffices to show that
$$\varepsilon_{z_i}(\partial_q^i(E_q^{\xi}(h_{n,q}(x)))=[n]_q!\delta_{i,n}.$$

Since two $q$-shift invariant operators commute, we have
$$\varepsilon_{z_i}(\partial_q^i(E_q^{\xi}(h_{n,q}(x)))=\varepsilon_{z_i}(E_q^{\xi}\partial_q^i(h_{n,q}(x)))
=\varepsilon_{z_i \oplus
\xi}(\partial_q^i(h_{n,q}(x)))=\varepsilon_{w_i}(\partial_q^i(h_{n,q}(x)))=[n]_q!\delta_{i,n}.\hfill\square$$
\end{dem}

\begin{prop}\label{shifttt}
Assume that $(h_{n,q}(x))_{n\geq0}$ be the generalized $q$-Gon$\breve{c}$arov basis
associated with the pair $(\partial_q,\mathcal{W})$, where
$\mathcal{W}=(w_i)_{i\geq0}$ is the grid given by $w_i=z_i \oplus (i \otimes \xi)$ with
$\xi$ is a fixed scalar and $i \otimes \xi=\underbrace{\xi\oplus...\oplus\xi}_{i\ terms}$.

 Then, $(h_{n,q}(x))_{n\geq0}$ is also the generalized $q$-Gon$\breve{c}$arov basis
 associated with the pair $(E_q^{\xi}\partial_q,\mathcal{Z})$.
\end{prop}
\begin{dem}
Let us notice that $(E_q^{\xi}\partial_q)^i=(E_q^{\xi})^i\partial_q^i$ and
$E_q^{a}E_q^{b}=E_q^{a\oplus b}$. Therefore
$$\varepsilon_{z_i}((E_q^{\xi}\partial_q)^i(h_{n,q}(x)))=\varepsilon_{z_i}E_q^{i\otimes\xi}(\partial_q^i(h_{n,q}(x)))=\varepsilon_{w_i}(\partial_q^i(h_{n,q}(x)))=[n]_q!\delta_{i,n}.\hfill\square$$
\end{dem}
For example, the q-analogue of the shifted factorial is
$$(x)_{(n,h,q)}=x(x-h[1]_q)...(x-h[n-1]_q).$$
This sequence is the basic sequence of the $q$-delta operator $\Delta_{h,q}$, given in \cite{7},
$$\Delta_{h,q}f(x)=\frac{f(qx+h)-f(x)}{(q-1)x+h}.$$
The generalized $q$-Gon$\breve{c}$arov polynomials associated with the $q$-delta $\partial_q=\Delta_{h,q}$ and
the $q$-arithmetic progression sequence $\mathcal{Z}$ = ($z_i=a\oplus i\otimes b$)$_{i\geq0}$
are the $\partial_q$-Abel polynomials of the form $t_{n,q}(x;\partial_q,(a\oplus i\otimes b)_{i\geq0})=(x\ominus a)\displaystyle{\prod_{i=1}^{n-1}(x\ominus a\ominus n\otimes b-h[i]_q ).}$\\
 In particular, the case $h=0$ gives the $q$-shift of the classical Abel polynomials associated to the sequence $(x^n)_{n\geq 0}$ and the same grid, for more details see \cite{4}.\\

The next proposition allows us to determinate a explicit formulas for the generalized
$q$-Gon$\breve{c}$arov basis from the explicit formulas of the basic sequence of the
$q$-delta operator $\partial_q$.
\begin{prop}\label{hhh}
Let $(t_{n,q}(x))_{n\geq0}$ be the generalized $q$-Gon$\breve{c}$arov basis associated with the
pair $(\partial_q,\mathcal{Z})$ and let $(p_n(x))_{n\geq0}$ be the sequence basic polynomials of the $q$-delta operator $\partial_q$. Then, for all $n \in \N$, we have
\begin{equation}\label{d(t)}
p_n(x)=\sum_{i=0}^n {n \brack i}_q p_{n-i}(z_i)t_{i,q}(x).
\end{equation}
Therefore
\begin{equation}\label{ttt}
t_{n,q}(x)=p_n(x)-\sum_{i=0}^{n-1} {n \brack i}_q p_{n-i}(z_i)t_{i,q}(x).
\end{equation}
\end{prop}
\begin{dem}
Applying (\ref{base}) and the definition of $\partial_q$ $i$ times on $p_n$($x$), we have
$$p_n(x)=\sum_{i=0}^{n}\frac{\varepsilon_{z_i}(\partial_q^i(p_n))}{[i]_q!}
t_{i,q}(x)=\sum_{i=0}^{n}\frac{(q;q)_n\ p_{n-i}(z_i)}{(q;q)_{n-i}\ (1-q)^i\ [i]_q!}
t_{i,q}(x) =\sum_{i=0}^n {n \brack i}_q p_{n-i}(z_i)t_{i,q}(x).$$
Thus, we obtain the desired equation (\ref{d(t)}).\hfill$\square$
\end{dem}

In the following result, we give the $q$-analogue of the generalized binomial expansion for
the generalized Gon$\breve{c}$arov polynomials.
\begin{prop}\label{identity}
Let $(t_{n,q}^{(j)}(x))_{n\geq0}$ be the generalized $q$-Gon$\breve{c}$arov basis associated
with the pair $(\partial_q,\mathcal{Z}^{(j)})$ and let $(p_n(x))_{n\geq0}$ be the sequence of basic polynomials of the $q$-delta operator $\partial_q$. Then, for all $\xi \in \mathbb{K}$ and $\ n\in \mathbb{N}$, we have
\begin{equation}\label{binomial}
t_{n,q}(x \oplus \xi)= \sum_{i=0}^n {n \brack i}_q\ t_{n-i,q}^{(i)}(\xi)\ p_i(x).
\end{equation}
In particular, if $\xi=0$, we have 
\begin{equation}\label{bino}
t_{n,q}(x)= \sum_{i=0}^n {n \brack i}_q\ t_{n-i,q}^{(i)}(0)\ p_i(x).
\end{equation}

\end{prop}
\begin{dem}
Let $\xi \in \mathbb{K}$ and $n\in \N$. Since $(p_i(x))_{0\leq i\leq n}$ is a basis of the linear subspace of $\mathbb{K}[x]$ of
polynomials of degree $\leq n$,  there exists $\{c_i, \;\ 0\leq i\leq n\}\subset \mathbb{K}$ such that
$$E_q^\xi (t_{n,q}(x))=\sum_{i=0}^n c_i p_i(x).$$
Thus, for any $0\leq j \leq n$, $$\varepsilon_{\xi}(\partial_q^jt_{n,q}(x)) = \varepsilon_{0}(\partial_q^jt_{n,q}(x \oplus \xi)) = \partial_q^j(E_q^\xi (t_{n,q}(x)))=\sum_{i=0}^n c_i\varepsilon_0(\partial_q^j p_i(x)).$$
The Eq. (\ref{defg}), $\partial_q^j t_{n,q}(x)=\frac{(q^{n-j+1},q)_{j}}{(1-q)^j}
t_{n-j,q}^{(j)}(x)$, implies that $\varepsilon_{\xi}(\partial_q^j t_{n,q}(x))=\frac{(q^{n-j+1},q)_{j}}{(1-q)^j}
t_{n-j,q}^{(j)}(\xi)$. 
Furthermore, the sequence $(p_n(x))_{n\geq0}$ being the  basic polynomials
of $\partial_q$, satisfies \\$\varepsilon_{0}(\partial_q^jp_i(x))=[i]_q\delta_{i,j}$, for all $i \in \N$.
Substituting all of the above, we obtain
$$c_j=\frac{(q^{n-j+1},q)_{j}}{(1-q)^j\ [j]_q!} t_{n-j,q}^{(j)}(\xi)= {n \brack j}_q\ t_{n-j,q}^{(j)}(\xi)$$
which gives the desired identity (\ref{binomial}).\hfill$\square$
\end{dem}

\begin{cor}
Assume  $\mathcal{Z}$ is a constant grid, means $z_i=z_j$ for all $i,j$, then the generalized
$q$-Gon$\breve{c}$arov basis associated with the pair $(\partial_q,\mathcal{Z})$ satisfies
\begin{equation}
t_{n,q}(x \oplus \xi)=\sum_{i=0}^n {n \brack i}_q\ t_{n-i,q}(\xi)\ p_i(x),
\end{equation}
and
\begin{equation}\label{q1}
t_{n,q}(x)=\sum_{i=0}^n {n \brack i}_q\ t_{n-i,q}(0)\ p_i(x),\;\ for \;\ \xi=0.
\end{equation}
\end{cor} \quad
\section{ A Combinatorial Formula For generalized $q$-Gon$\breve{c}$arov polynomials}

 In this section, we shall give a $q$-analogue combinatorial formula of the constant terms (terms without $p_n$($x$), for all $n\geq 1$) of the generalized $q$-Gon$\breve{c}$arov polynomial associated with the pair ($\partial_q,\mathcal{Z}$). This interpretation
is a generalisation of the techniques used in \cite{1} and \cite{2}.\\
In view of Eq. (\ref{ttt}), the first few generalized $q$-Gon$\breve{c}$arov polynomials are:
\begin{align*}
t_{0,q}(x,\partial_q,\mathcal{Z})&=1,\\
t_{1,q}(x,\partial_q,\mathcal{Z})&=p_1(x)-p_1(z_0),\\
t_{2,q}(x,\partial_q,\mathcal{Z})&=p_2(x)-[2]_q p_1(z_1)p_1(x)+[2]_q p_1(z_0)
p_1(z_1)-p_2(z_0),\\
t_{3,q}(x,\partial_q,\mathcal{Z})&=p_3(x)-[3]_q p_1(z_2) p_2(x)+\big([3]_q[2]_q p_1(z_1)
p_1(z_2)-[3]_q p_2(z_1)\big)p_1(x)\\
&+[3]_q p_1(z_0) p_2(z_1)-[3]_q[2]_q p_1(z_2) p_1(z_1) p_1(z_0) + [3]_q p_2(z_0) p_1(z_2)
-p_3(z_0).
\end{align*}

It is well known that if $q$ tends to $1$, those results returns to the formulas of generalized
Gon$\breve{c}$arov polynomial $t_n(x,\partial,\mathcal{Z})$ obtained by Lorentz, Tringali and Yan in their paper \cite{2}.\\

From an ordered partition $\rho=(B_1,...,B_k)$ of a set with $n$ elements, we set
$|\rho|=k$ and we define, for every $i=1,...,k$, $b_i$ to be the size of the block $B_i$
and $s_i=\sum_{j=1}^i b_{j}$ with $s_0=0$.

 For a fixed $n$, we will choose among the ordered partitions of $n$, the partitions whose
 elements of their blocks are ordered from $1$ to $n$. For example, we take among the
 ordered partitions of $3$, the following partitions
$(\{1\},\{2\},\{3\})$, $(\{1,2\},\{3\})$, $(\{1\},\{2,3\})$ and $(\{1,2,3\})$.
We denote by [$n$] the set $\{1, 2, \ldots, n\}$, $\mathfrak{\mathfrak{P}}_n$ the set of
the all ordered partitions under hypothesis above of [$n$] and $\mathfrak{\mathfrak{P}}_n^k$
the subset of $\mathfrak{\mathfrak{P}}_n$ such that all partition has $k$ blocs.

Under these hypothesis, we  give in the next theorem a combinatorial formula of the
constant terms of generalized $q$-Gon$\breve{c}$arov polynomials associated to the pair
($\partial_q$, $\mathcal{Z}$), which is a $q-$analog of Theorem 4.1 in
\cite{2}.
\begin{thm}
Let $(t_{n,q}(x))_{n\geq0}$ be the generalized $q$-Gon$\breve{c}$arov basis associated with
the pair $(\partial_q,\mathcal{Z})$ and $(p_n(x))_{n\geq0}$ be the sequence of basic
polynomials of $\partial_q$. Then, for any $n\geq1$, we have
\begin{equation}\label{equation1}
t_{n,q}(0)=\sum_{k=1}^n\sum_{\rho \in \mathfrak{\mathfrak{P}}_n^k} (-1)^k
\prod_{i=0}^{k-1}\  {s_{k-i} \brack b_{k-i}}_q\ p_{b_{k-i}}(z_{n-s_{k-i}}).
\end{equation}
\end{thm}
\begin{dem}
Using the Proposition \ref{hhh} and by considering that $p_n(0)=0$, for all $n\geq 1$, we
have
\begin{equation}\label{equation2}
t_{n,q}(0)=-\sum_{i=0}^{n-1} {n \brack i}_q p_{n-i}(z_i)\ t_{i,q}(0).
\end{equation}
We will reason by recurrence on $n$ and we denote the right-hand side of Eq.
(\ref{equation1}) by $\mathcal{T}_n(0)$, for $n\geq 1$. We let $\mathcal{T}_0(0)=1$ and
we suppose that for any $1\leq i \leq n-1$ and in order to the ordered partitions on [$i$], we
have
$$\mathcal{T}_{i}(0)=\sum_{k=1}^i\sum_{\rho \in \mathfrak{\mathfrak{P}}_i^k} (-1)^k
\prod_{j=0}^{k-1}\  {s_{k-j} \brack b_{k-j}}_q\ p_{b_{k-j}}(z_{n-s_{k-j}}).$$

For $n=1$, it's easily to see that $t_{1,q}(0)=\mathcal{T}_1(0)$. So suppose that the
equality is true until $n-1$ and we will show that it is true for $n$.\\

The idea of this proof comes from the fact that to prove $\mathcal{T}_n=t_{n,q}$, it
suffices to determinate the coefficient of $p_{n-i}(z_i)$ in $\mathcal{T}_n(0)$ which is
equivalent to those identities: $$b_{k-i}=n-i\ and\ n-s_{k-i}=i\Longleftrightarrow k-i=1\
and\ b_1=s_1=n-i.$$

 First, we have
 $\mathfrak{P}_n=\displaystyle{\bigcup_{k=1}^n\mathfrak{P}_n^k=\bigcup_{j=0}^{n-1}\mathfrak{P}_{n,j}}$,
 where $\mathfrak{P}_{n,j}=\{\rho\in \mathfrak{P}_n, \; |B_1|=n-j\}$.

On the other hand, for all $0 \leq i\leq n-1$, let $X : = \{a_1, a_2, \ldots, a_{n-i}\}$ a
subset of [$n$] such that $a_1 < a_2 < \ldots < a_{n-i}$. We defined the ordered partition
on the set $X$ like the one on the set [$n-i$]. Thus, we get a bijection between $\{
(X,\rho'), \rho' \in  \mathfrak{P}_i\}$ and $\mathfrak{P}_{n,i} $.

So, in the sense of $q$-analog, we have ${n \brack i}_q$ ways to choose a subset $X$ of
$n-i$ elements of [$n$]. The fact that the rest [$n$]$\backslash X$ is seen as well as the
set [$i$], and from the recurrence hypothesis on [$i$], we have the coefficient of
$p_{n-i}(z_i)$ in $\mathcal{T}_n(0)$ is
$$-{n \brack i}_q\sum_{k=1}^i\sum_{\rho \in \mathfrak{\mathfrak{P}}_i^k} (-1)^k
\prod_{j=0}^{k-1}\  {s_{k-j} \brack b_{k-j}}_q\ p_{b_{k-j}}(z_{n-s_{k-j}})= -{n \brack
i}_q\mathcal{T}_{i}(0).$$

By summing over $i = 0, 1, \ldots , n - 1 $,  we obtain that $\mathcal{T}_n$ verifies the
equation (\ref{equation2}). Whence, we have the desired recurrence.\hfill$\square$
\end{dem}

As a reinforcement of this demonstration, we will give a small example to better
understand. Let's put $n=4$. To prove that $\mathcal{T}_4(0)=t_4(0)$, it suffices to show
that $\mathcal{T}_4 (0)$ verifies Eq. (\ref{equation2}) for $n=4$. So, we must search the
coefficients of $p_{4-j}(z_j)$, for $0\leq j\leq 3$ in $\mathcal{T}_4 (0)$.  We have
\begin{align*}
\mathcal{T}_4(0)&=\big([4]_qp_3(z_0)-{4 \brack 2}_q[2]_q p_2(z_0)p_1(z_2)-[4]_q[3]_q
p_1(z_0)p_2(z_1)+[4]_q!p_1(z_0)p_1(z_1)p_1(z_2)\big)\underline{p_1(z_3)}\\
&+\big({4 \brack 2}_q p_2(z_0)-[4]_q[3]_q p_1(z_0)p_1(z_2)\big)\underline{p_2(z_2)}
+[4]_q p_1(z_0)\underline{p_3(z_1)} -\underline{p_4(z_0)}.
\end{align*}
 First, we take $j=3$ and we will choose one element of $[4]$. By simple calculation, we find
 that the coefficient of $p_{4-3}(z_3)=p_1(z_3)$ is $-{4 \brack 3}_q  \mathcal{T}_3(0)$.
 Similarly, if we take $j=2$ then we obtain that the coefficient of $p_2(z_2)$ is $-{4 \brack
 2}_q \mathcal{T}_2(0)$.
We will determine the calculation for $j=1$ and $j=0$ until finally finds that
$$\mathcal{T}_n(0)=-\sum_{j=0}^{3}{4 \brack j}_q p_{4-j}(z_j)\mathcal{T}_j(0).$$
\quad

Now, we focus on the constant terms of generalized $q$-Gon$\breve{c}$arov polynomials  and we count the number  $f_{n,q}$ of monomials in this
constant term. The first few values of  $f_{n,q}$ are\\

$f_{0,q}= 1$,

$f_{1,q}= 1$,

$f_{2,q}= 2+q$,

$f_{3,q}=1+[3]_q(3+q)$,

$f_{4,q}=1+[4]_q(2+[3]_q!+2[3]_q)+(1+[2]_q){4 \brack 2}_q $,

$f_{5,q}=1+[5]_q$\big(2+[4]$_q$!+2[4]$_q$+2[4]$_q$[3]$_q$+${4 \brack 2}_q$(2+q)\big)+${5 \brack
  2}_q$(2+q)+${5 \brack 3}_q$\big(1+2[3]$_q$+[3]$_q$!\big).
\quad\\

From this table, we can find, in the next proposition, the number $f_{n,q}$ which can be
considered as a $q$-analog of Fubini's number $f_n$, refer to \cite{8}.
\begin{prop} The number $f_{n,q}$ satisfy the following recurrence:
\begin{equation}\label{equation7}
f_{n,q}=\sum_{k=0}^{n-1}{n \brack k}_q f_{k,q},
\end{equation}
where $f_{0,q}=1$.
\end{prop}
\begin{dem} Substituting $\mathcal{T}_i$(0) by $f_{i,q}$ and replacing $p_{n-i}(z_i)$ by $-1$, for all $0\leq i \leq n-1$ and for any $z_i$, in the above proof, we obtain the desired equation
(\ref{equation7}).\hfill$\square$
\end{dem}
\quad\\

Hence, we obtain the following result which expresses a new formula of $q$-Fubini numbers in terms of order partitions.
\begin{prop} Under the hypothesise on the order partitions above and for $n \geq 0$, we have
\begin{equation}
f_{n,q}=\sum_{k=1}^n\sum_{\rho \in \mathfrak{\mathfrak{P}}_n^k} \prod_{i=0}^{k-1}\
{s_{k-i} \brack b_{k-i}}_q,
\end{equation}
where $f_{0,q}=1$.
\end{prop}

\bibliographystyle{amsplain}


\end{document}